\documentclass[12pt]{amsart}
\usepackage[spanish]{babel}
\usepackage[spanish]{babelbib}
\usepackage{mdwlist}
\usepackage[utf8]{inputenc}

\begin{document}

\title[Sobre la ecología social...]{Sobre la ecología social de las instituciones de educación superior en la Cañada oaxaqueña}
\author{Octavio Alberto Agustín-Aquino}
\address{Universidad de la Ca\~nada\\
San Antonio Nanahuatipan Km 1.7 s/n. Paraje Titlacuatitla\\
Teotitlan de Flores Mag\'on, Oaxaca, M\'exico, C.P. 68540.}
\email{octavioalberto@unca.edu.mx}
\begin{abstract}
Se proporcionan algunas observaciones en relación a la competencia que hay entre otras universidades y
la Universidad de la Cañada (Unca) usando el modelo de Lotka-Volterra respecto a la carrera de Informática, empleando las ideas introducidas por
McPherson en 1983. Como consecuencia, se señalan algunas líneas de acción a seguir para la Unca.
\end{abstract}
\keywords{Matrícula, competencia, universidades, modelo de Lotka-Volterra}
\renewcommand{\datename}{\textit{Fecha}:}
\renewcommand{\emailaddrname}{\textit{Correo electrónico}}
\renewcommand{\keywordsname}{\textit{Palabras y frases clave}}
\subjclass[2010]{92D40, 97B40}
\date{5 de junio de 2015}
\maketitle

\section{Introducción}
La Universidad de la Cañada (Unca) es una integrante del Sistema de Universidades del Estado de Oaxaca (Suneo),
y se localiza en Teotitlán de Flores Magón, en la región de la Cañada. Pese a que en la región es la que tiene el
mayor rezago educativo y la menor participación económica en el estado \cite[Sec. 8.2.1]{ped1116}, la matrícula proveniente de las poblaciones circunvecinas a la Unca parece ser notablemente baja
en el área de Informática, considerando que la institución es estatal, prácticamente gratuita y que la demanda laboral en esta área
es bastante grande.

Sin embargo, resulta que en términos del modelo de Lotka-Volterra para la competencia entre dos o más especies, puede
encontrarse parte de una explicación a este fenómeno, y además permite dar algunas pautas para adoptar políticas más eficaces para incrementar
su alumnado. También nos permite formarnos un panorama para el futuro una vez que se abra un campus de Nova Universitas,
otra universidad perteneciente al Suneo, en Huautla de Jiménez.

\section{El modelo para dos competidores}

Miller McPherson introdujo en \cite{mM83} la idea de aplicar a las sociedades de afiliación voluntaria el mismo modelo que Alfred Lotka y Vito
Volterra desarrollaran \cite{SZ78} para entender la competencia entre organismos de especies diferentes por un mismo recurso, y parece
pertinente hacer lo mismo para analizar la dinámica de la matrícula universitaria cuando una institución compite contra otra. En el caso
de la Unca, puesto que es la única institución de educación superior en la región de la Cañada en Oaxaca, tiene sentido abstraer
al resto de las instituciones que compiten con ella por estudiantes como una sola entidad.

Las ecuaciones de Lotka-Volterra para dos especies, asociaciones o universidades que compiten son
\begin{align*}
\frac{d}{dt} N_{1} &= r_{1}N_{1}\left(1-\frac{N_{1}+\alpha_{1,2}N_{2}}{K_{1}}\right),\\
\frac{d}{dt} N_{2} &= r_{2}N_{2}\left(1-\frac{N_{2}+\alpha_{2,1}N_{1}}{K_{2}}\right),
\end{align*}
donde $N_{i}$ es la población, $K_{i}$ es la capacidad de carga y $r_{i}$ es la
tasa de reproducción de la especie $i$-ésima\footnote{Para comprender mejor los significados de estos términos, puede
consultarse \cite[Capítulo 3]{lR15} y el artículo de McPherson \cite{mM83}.}. Los coeficientes $\alpha_{i,j}$ representan la competencia
que hay entre ambas especies. En términos biológicos, puede pensarse como cuántos elementos de la especie $j$-ésima
pueden canjearse por uno de la especie $i$-ésima respecto a los mismos recursos.

Para un sistema de dos competidores, hay cuatro tipos de equilibrio \cite[Cap. 5]{BCC12}
\begin{enumerate*}
\item En dos de ell, los únicos equilibrios aparte del trivial son que sólo una de las especies se extingue.
\item En otro, existe un equilibrio adicional donde ambas coexisten, pero es intestable.
\item En un cuarto, ambas especies coexisten y el equilibrio es estable.
\end{enumerate*}

Según lo observado en la Unca, el número de alumnos provenientes del municipio de Teotitlán no ha cambiado apreciablemente
en los últimos cinco años\footnote{Esta es una observación directa del autor.}, lo que permite suponer que el sistema tiene
un equilibrio estable. Por lo tanto, es viable aproximar las poblaciones $N_{i}^{*}$ en equilibrio y las capacidades de carga
$K_{i}$ a partir de un estimado de los coeficientes $a_{i,j}$.

Para realizar este tipo de estimación, McPherson determina los llamados \emph{nichos} de las organizaciones
en términos de las características de edad y nivel educativo, que son llamadas \emph{dimensiones} del nicho. Para el caso de las universidades, el nicho en cuanto a edad y
nivel educativo es sumamente angosto e idéntico entre los competidores (los alumnos de bachillerato mayores de edad). Por lo tanto, es necesario considerar una dimensión
alternativa y más decisiva, que es el nivel socioeconómico de las familias de los estudiantes. Tomaremos como referencia el porcentaje de las familias que tienen un ingreso mayor o igual
a tres salarios minímos, que es alrededor de una cuarta parte\footnote{El porcentaje de los que ganan más de tres salarios
mínimos respecto al total que reporta algún ingreso es $21{.}31\%$. Si se restringe a la fracción que gana más de cinco
salarios mínimos respecto a los que ganan más de tres, el porcentaje es $25{.}60\%$.} respecto al total que reporta algún ingreso en Oaxaca \cite[Cuadro 8.4]{AEGPEF13}. La hipótesis de trabajo es que las familias con suficiente ingreso pueden decidir enviar a sus hijos a estudiar
fuera, mientras que el resto tiene un fuerte incentivo para inscribirlos en la universidad local.

Siguiendo el esquema de McPherson, esto significa que $\alpha_{1,2} = 0{.}25$ y $\alpha_{2,1} = 1$, donde tomamos
al índice $2$ para la Unca y $1$ para sus competidores. Estos números, en términos sociológicos, significan que los competidores externos como pueden perder toda
su matrícula potencial de Teotitlán con la Unca y que ésta puede perder hasta una cuarta parte de la suya contra la competencia.
Por otro lado, de la matrícula de la Unca se observa que su población en equilibrio en la carrera de Informática es $N_{2}^{*}= 8$, y
se estima que $N_{1}^{*} = 24$, considerando los datos del Inegi para la matrícula universitaria de Oaxaca, Puebla y Veracruz \cite[Cuadro 4.19]{AEGPEF13}.

Con estos datos, igualando a cero las razones instantáneas de crecimiento de las poblaciones estudiantiles y resolviendo el sistema,
se encuentra que
\begin{align*}
K_{1} &= N_{1}^{*}+\alpha_{1,2}N_{2}^{*} = 24+\tfrac{1}{4}\times 8 = 26,\\
K_{2} &= N_{2}^{*}+\alpha_{2,1}N_{1}^{*} = 8+1\times 24 = 32.
\end{align*}

Esto significa que la Unca está, al menos en teoría\footnote{Debe señalarse que McPherson encontró en su estudio que el modelo así
planteado normalmente exagera la capacidad de carga. Pese a ello, predice correctamente las magnitudes relativas de los nichos.
Con esto en mente debe observarse, en la circunstancia bajo estudio, que las capacidades predichas indican que la Unca está ligeramente
por encima de sus competidores en este aspecto. Aunque este hecho hace más paradójicas las preferencias de los estudiantes para continuar
sus estudios en otros lugares, también es otro indicador de la calidad de la Unca.} y disponiendo de una
situación favorable, en posibilidades de multiplicar significativamente
su matrícula local (esto es, la proveniente del municipio donde se encuentra). Sin embargo, esto presenta un indicio de
que la competencia con el exterior es algo intensa y que la Unca no está convenientemente
preparada para asumir este reto\footnote{Hay que considerar que la gran mayoría de los estudiantes de la Unca de la carrera de Informática proviene de municipios distintos al de Teotitlán;
si se sumaran a la matrícula local obtenida aquí, con seguridad alcanzarían el cupo de dos aulas. Hay que añadir que el Departamento de Servicios Escolares de la Unca desde hace algún
tiempo enfrenta notables dificultades para distribuir los espacios y tiempos, evidenciadas por la introducción de un programa que calcula
la configuración óptima de los horarios en 2013, según me informó su titular.}.

\section{El modelo para tres competidores}

Otra situación\footnote{Agradezco a mi colega María del Rosario Peralta Calvo el atraer mi atención hacia este tema.} que
potencialmente puede ser problemática para la Unca es la apertura del campus de Nova Universitas
en Huautla de Jiménez, ubicada a unos 26 km de Teotitlán en línea recta. Aunque también pertenece
al Suneo, representa una opción muy atractiva para las poblaciones más pobres que circundan a Huautla.

Las ecuaciones de Lotka-Volterra para este caso son
\begin{align*}
\frac{d}{dt} N_{1} &= r_{1}N_{1}\left(1-\frac{N_{1}+\alpha_{1,2}N_{2}+\alpha_{1,3}N_{3}}{K_{1}}\right),\\
\frac{d}{dt} N_{2} &= r_{2}N_{2}\left(1-\frac{N_{2}+\alpha_{2,1}N_{1}+\alpha_{2,3}N_{3}}{K_{2}}\right),\\
\frac{d}{dt} N_{3} &= r_{3}N_{3}\left(1-\frac{N_{3}+\alpha_{3,1}N_{1}+\alpha_{3,2}N_{2}}{K_{3}}\right),
\end{align*}
con los mismos significados de los parámetros, y tomando el índice 3 para Nova Universitas. Al introducir una nueva
variable el sistema tiene, en lugar de cuatro, ocho equilibrios posibles, y todos ellos pueden presentarse con una elección
apropiada de los parámetros. Además, es significativamente más difícil determinar sus características (véase \cite{ML75}
para un ejemplo concreto y \cite{sS76} para los detalles finos).

Sin embargo, no renunciamos a obtener información del modelo, si bien puede ser muy tentativa. Por principio de cuentas,
no es descabellado suponer que los coeficientes de competencia del campus de Nova contra el resto de universidades fuera de la región
son los mismos que los de la Unca, es decir, $\alpha_{3,1} = \alpha_{2,1}$ y $\alpha_{1,3}=\alpha_{1,2}$, pues las
características de la población de la Cañada en total se mantienen.
Sin embargo, un problema con la estimación de los coeficientes restantes es que la dimensión socioeconómica no determina
a $\alpha_{2,3}$ y $\alpha_{3,2}$, pues ambas universidades son estatales, públicas y de calidad semejante.

En tal caso, lo más razonable es considerar que la zona mazateca donde se encuentra Huautla de Jiménez tiene una
población bastante mayor que la de las cercanías de Teotitlán. Parece razonable, pues, estimar que
\[
 \alpha_{3,2} = \frac{8966}{30004} \approx 0{.}3
\]
es el cociente de las poblaciones respectivas \cite{cpv10}, y también que $\alpha_{2,3}=1-\alpha_{2,3}= 0{.}7$. Suponiendo
ahora que las capacidades de carga estimadas en la sección anterior se mantienen y que la del campus de Nova Universitas
será de alrededor de la mitad del de la Unca\footnote{Esto se debe a que el modelo de Nova Universitas
es semipresencial, donde los alumnos del campus periférico interactúan por video-conferencia con el profesor ubicado
en el campus central de Ocotlán de Morelos.}, realizamos una simulación numérica de
$1500$ iteraciones con tamaño de paso $0{.}01$, tomando las poblaciones iniciales $N_{1}=24$, $N_{2}=8$ y $N_{3}=24$
(un panorama muy optimista) e idénticas tasas intrínsecas de crecimiento $r_{i}=1$. El resultado es que el sistema evoluciona a un equilibrio
donde las poblaciones de la Unca y sus competidores fuera de la Cañada
regresan a su estado original, mientras que la del campus de Nova desaparece.

Ajustando por medio de más simulaciones los parámetros, parece que cuando la capacidad de carga del campus de Nova
es de alrededor de $K_{3}=29$ puede sobrevivir durante un tiempo razonable. Además, si $K_{3}=31$ (prácticamente la misma capacidad estimada para
la Unca), parece salir a flote sin reducir significativamente la matrícula de la Unca y, en todo caso, aminora la salida de alumnos
de la región\footnote{Una simulación de quince mil iteraciones con las mismas condiciones iniciales y tamaño de paso ya mencionadas con
tal capacidad de carga, arroja el resultado $N_{1}^{*}=23{.}5$, $N_{2}^{*}=7{.}9$ y $N_{3}^{*}=1{.}9$.}.

\section{Estrategias y conclusiones}

Los parámetros que puede controlar la Unca, suponiendo que la situación de los competidores se mantiene igual, son la capacidad de carga $K_{2}$ y el coeficiente
de competencia $\alpha_{1,2}$. El primero implica la ampliación de las instalaciones (aulas, laboratorios, salas de cómputo, cafetería) y de la planta laboral
(profesores y personal administrativo) para atender a un número mayor de alumnos, mientras que el segundo se relaciona con aumentar el atractivo que representa la universidad respecto a las alternativas (a través de más becas, mayor difusión y vinculación). Sin embargo,
al tomar las derivadas de la población $N_{2}^{*}$ en equilibrio respecto a estos parámetros,
se descubre que la razón de cambio es constante al cambiar a $K_{2}$, mientras que depende del valor actual de $\alpha_{2,1}$ para $\alpha_{2,1}$ mismo.
Cuando se consideran incrementos pequeños y solamente uno de estos parámetros puede cambiarse, resulta que  el valor de $N_{2}^{*}$ crece más aumentando
la capacidad de carga que disminuyendo el coeficiente de competitividad. Por ser más específicos: si $K_{2}$ aumenta en $1\%$, entonces $N_{2}^{*}  \approx 8{.}42$,
mientras que si $\alpha_{2,1}$ disminuye (y por ello mejora, pues significa que un menor porcentaje de alumnos decide estudiar fuera) en la misma proporción, se tiene $N_{2}^{*}\approx 8{.}29$.
Dicho en otras palabras: para incrementar en por lo menos un alumno el equilibrio, se debe incrementar en $2{.}4\%$ la capacidad de carga o disminuir en $3{.}5\%$
el coeficiente de competitividad.

Si bien se puede argumentar que el valor del incremento es tan pequeño y el modelo tan aproximado\footnote{En particular, es continuo.}
que no se puede distinguir con la suficiente precisión la mejor
estrategia, también es difícil no reconocer a partir del mismo que los cambios en las condiciones de los competidores, aún si son pequeños,
tiene consecuencias serias para la Unca. No es inverosímil esperar que la capacidad de recibir alumnos en Tehuacán, por citar un caso concreto pero significativo, crezca
en un $10\%$ en un año determinado, lo que ocasionaría una reducción de la matrícula en equilibrio a casi la mitad de alumnos para la Unca. Esto favorece la conclusión de que es vital
para el futuro de la universidad que amplíe su infraestructura para enfrentar mejor los cambios en su entorno.

Por otra parte, es comprensible que no sea económicamente factible llevar a cabo las mejoras en el rubro que más impacto tiene, y por lo tanto
sea necesario concentrarse en promover a la universidad para atraer a más estudiantes. Del modelo se colige, en esta dirección, que resulta
particularmente ventajoso dirigirse a las comunidades con mayor grado de marginación, donde la Unca representa una oportunidad real para acceder
a la educación superior para sus jóvenes, siendo que el equilibrio en Teotitlán ya fue alcanzado.

Finalmente, pese a las complejidades inherentes al modelo con tres competidores, parece sensato concluir de un análisis 
preliminar que la apertura de un campus de Nova Universitas no representa una amenaza para la Unca; de hecho, si su capacidad
de carga no es lo suficientemente grande, estaría destinada a la extinción por lo menos en la carrera de Informática.
Por otro lado, si se le proyectase como una universidad de envergadura comparable a la de la Unca, entonces
no sólo sobrevive, sino que incluso ayuda a acrecentar la población estudiantil local de la Unca y a impedir la descapitalización social
de la Cañada sobre la que hace énfasis el rector del Suneo.

\bibliographystyle{amsplain}
\bibliography{reporte19}

\providecommand{\bysame}{\leavevmode\hbox to3em{\hrulefill}\thinspace}
\providecommand{\MR}{\relax\ifhmode\unskip\space\fi MR }
\providecommand{\MRhref}[2]{%
  \href{http://www.ams.org/mathscinet-getitem?mr=#1}{#2}
}
\providecommand{\href}[2]{#2}
\begin{thebibliography}{1}

\bibitem{cpv10}
\emph{Censo de poblaci\'{o}n y vivienda}, INEGI, 2010.

\bibitem{ped1116}
\emph{Plan estatal de desarrollo de {O}axaca 2011-2016}, Gobierno
  Constitucional del Estado Libre y Soberano de Oaxaca, 2011.

\bibitem{AEGPEF13}
\emph{Anuario estad\'{\i}stico y geogr\'{a}fico por entidad federativa}, INEGI,
  2013.

\bibitem{BCC12}
Fred Brauer and Carlos Castillo-Ch\'{a}vez, \emph{Mathematical models in
  population biology and epidemiology}, 2nd ed., Texts in Applied Mathematics,
  vol.~40, Springer New York, 2012.

\bibitem{ML75}
Robert~M. May and Warren~J. Leonard, \emph{Nonlinear aspects of competition
  between three species}, SIAM Journal on Applied Mathematics \textbf{29}
  (1975), no.~2, 243--253.

\bibitem{mM83}
Miller McPherson, \emph{An ecology of affiliation}, American Sociological
  Review \textbf{48} (1983), no.~4, 519--532.

\bibitem{lR15}
Larry~R. Rockwood, \emph{Introduction to population ecology}, 2nd ed., John
  Wiley \& Sons, 2015.

\bibitem{SZ78}
Francesco~M. Scudo and James~R. Ziegler (eds.), \emph{The golden age of
  theoretical ecology: 1923-1940}, Lecture Notes in Biomathematics, vol.~22,
  Springer Berlin Heidelberg, 1978.

\bibitem{sS76}
Stephen Smale, \emph{On the differential equations of species in competition},
  Journal of Mathematical Biology \textbf{3} (1976), no.~1, 5--7.

\end{thebibliography}

\end{document}